\documentclass{gtart_h}

\def\ifplaintex{\expandafter\ifx\csname documentclass\endcsname\relax}

\def\gtp{{\mathsurround=0pt\it $\cal G\mskip-2mu$eometry \&\ 
$\cal T\!\!$opology $\cal P\!$ublications}}  

\def\Addressesr{\bigskip
{\small \parskip 0pt \leftskip 0pt \rightskip 0pt plus 1fil \def\\{\par}
\sl\theaddress\par
\medskip
\rm Email:\stdspace\tt\theemail\hfill\rm Received:\qua\receiveddate \par}}

\def\recd{{\small Received:\qua\receiveddate\ifx\reviseddate\relax
\else\qquad Revised:\qua\reviseddate\fi\par}} 

\def\lognumber#1{\def\thelognumber{#1}}
\def\volumenumber#1{\def\thevolumenumber{#1}}
\def\volumeyear#1{\def\thevolumeyear{#1}}
\def\papernumber#1{\def\thepapernumber{#1}}
\def\pagenumbers#1#2{\def\startpage{#1}\def\finishpage{#2}}
\def\published#1{\def\publishdate{#1}}

\def\received#1{\def\receiveddate{#1}}

\def\accepted#1{\def\accepteddate{#1}}

\def\asciiaddress#1{\def\theasciiaddress{#1}}

\let\\\par\let\thelognumber\relax\let\thevolumenumber\relax
\let\thepapernumber\relax\let\thevolumeyear\relax\let\startpage\relax
\let\finishpage\relax\let\publishdate\relax\let\receiveddate\relax
\let\reviseddate\relax\let\accepteddate\relax\let\theasciititle\relax
\let\theasciiauthors\relax\let\theasciiaddress\relax
\let\theasciiabstract\relax

\let\theasciiemail\relax

\ifplaintex
\font\logobig=cmssbx10 scaled 3836
\font\logomed=cmssbx10 scaled 2557
\else
\font\logobig=cmssbx10 scaled 4200
\font\logomed=cmssbx10 scaled 2800
\fi

\long\def\makeagttitle{   
\count0=\startpage
\agt\hfill      
\hbox to 45truept{\vbox to 0pt{\vglue -13truept{\logomed A\kern -.37em{\logobig 
T}\kern -.38em G}\vss}\hss}
\break
{\small Volume \thevolumenumber\ (\thevolumeyear)
\startpage--\finishpage\nl
Published: \publishdate}

\vglue .25truein

{\parskip=0pt\leftskip 0pt plus
1fil\def\\{\par\smallskip}{\Large\bf\thetitle}\par\medskip} \vglue
0.05truein

{\parskip=0pt\leftskip 0pt plus 1fil\def\\{\par}{\sc\theauthors}
\par\medskip}%
 
\vglue 0.03truein 

{\small\leftskip 25truept\rightskip 25truept{\bf Abstract}\stdspace\theabstract

{\bf AMS Classification}\stdspace\theprimaryclass
\ifx\thesecondaryclass\relax\else; \thesecondaryclass\fi\par
{\bf Keywords}\stdspace \thekeywords\par}\vglue 7truept

}   

\ifplaintex
\font\phead=cmsl9 scaled 950
\font\pnum=cmbx10 scaled 913
\font\pfoot=cmsl9 scaled 950
\headline{\vbox to 0pt{\vskip -4.5mm\line{\small\phead\ifnum
\count0=\startpage ISSN 1472-2739 (on-line) 1472-2747 (printed)
\hfill {\pnum\folio}\else\ifodd\count0\def\\{ }%
\ifx\theshorttitle\relax\thetitle\else\theshorttitle\fi\hfill{\pnum\folio}
\else\def\\{ and }{\pnum\folio}\hfill\ifx\theshortauthors\relax\theauthors
\else\theshortauthors\fi\fi\fi}\vss}}
\footline{\vbox to 0pt{\vglue 0mm\line{\small\pfoot\ifnum\count0=\startpage
\copyright\ \gtp\hfill\else
\agt, Volume \thevolumenumber\ (\thevolumeyear)\hfill\fi}\vss}}
\else
\font=cmsl9 scaled 1050
\font\lnum=cmbx10 
\font=cmsl9 scaled 1050
\makeatletter
\def\@oddhead{{\small\lhead\ifnum\count0=\startpage ISSN 1472-2739 
(on-line) 1472-2747 (printed)\hfill {\lnum\number\count0}\else\ifodd\count0
\def\\{ }\ifx\theshorttitle\relax \thetitle \else\theshorttitle\fi\hfill
{\lnum\number\count0}\else\def\\{ and }{\lnum\number\count0}
\hfill\ifx\theshortauthors\relax 
\theauthors\else\theshortauthors\fi\fi\fi}}\def\@evenhead{\@oddhead}
\def\@oddfoot{\small\lfoot\ifnum\count0=\startpage\copyright\ \gtp\hfill\else
\agt, Volume \thevolumenumber\ (\thevolumeyear)\hfill\fi}
\def\@evenfoot{\@oddfoot}
\makeatother
\fi
\let\maketitlepage\makeagttitle
\let\makeshorttitle\maketitlepage
\let\maketitle\maketitlepage

\newwrite\gtoutfile
\long\gdef\makeheadfile{  
{\def\\{, }\def\s{ }
\immediate\openout\gtoutfile head.xxx
\immediate\write\gtoutfile{Proxy-for: \ifx\theasciiauthors\relax
\theauthors\else\theasciiauthors\fi\s<\ifx\theasciiemail\relax\theemail\else\theasciiemail\fi>}
\immediate\write\gtoutfile{\noexpand\\}
\immediate\write\gtoutfile{Authors: \ifx\theasciiauthors\relax
\theauthors\else\theasciiauthors\fi}
{\def\\{ }\immediate\write\gtoutfile{Title: \ifx\theasciititle\relax
\thetitle\else\theasciititle\fi}}
\immediate\write\gtoutfile{Subj-class: GT or SG, GR etc}
\immediate\write\gtoutfile{MSC-class: \theprimaryclass\ifx\thesecondaryclass\relax\else, \thesecondaryclass\fi}
\immediate\write\gtoutfile{Journal-ref: Algebr. Geom. Topol. \thevolumenumber\s
(\thevolumeyear) \startpage-\finishpage}
\immediate\write\gtoutfile{Comments: Published by Algebraic and
Geometric Topology at}
\immediate\write\gtoutfile{\s\s\s  http://www.maths.warwick.ac.uk/agt/AGTVol\thevolumenumber/agt-\thevolumenumber-\thepapernumber.abs.html}
\immediate\write\gtoutfile{\noexpand\\}
\immediate\write\gtoutfile{}
\ifx\theasciiabstract\relax
\immediate\write\gtoutfile{\theabstract}\else
\immediate\write\gtoutfile{\theasciiabstract}\fi
\immediate\write\gtoutfile{}
\immediate\write\gtoutfile{\noexpand\\}
\immediate\write\gtoutfile{}
\immediate\closeout\gtoutfile}}  

\def\maketitlepage{\makeagttitle\makeheadfile}
\let\makeshorttitle\maketitlepage
\let\maketitle\maketitlepage

\lognumber{4}
\volumenumber{5}
\volumeyear{2005}
\papernumber{4}
\pagenumbers{53}{69}
\received{6 December 2004} 
\accepted{7 January 2005}
\published{9 January 2005}

\usepackage{amsmath,amssymb}
\usepackage[all]{xy}

\newcommand{\ZZ}{\mathbb{Z}}
\newcommand{\QQ}{\mathbb{Q}}

\newcommand{\PP}{\mathcal{P}}
\newcommand{\legendre}{\mathcal{L}}

\newcommand{\sym}{\mathfrak{S}}

\newcommand{\sous}{\curvearrowleft}
\newcommand{\id}{\operatorname{Id}}
\newcommand{\ch}{\operatorname{ch}}
\newcommand{\Ch}{\operatorname{\mathbf{Ch}}}

\newcommand{\Det}{\operatorname{Det}}
\newcommand{\comm}{\operatorname{Comm}}
\newcommand{\asso}{\operatorname{Assoc}}
\newcommand{\lie}{\operatorname{Lie}}
\newcommand{\dias}{\operatorname{Dias}}
\newcommand{\dend}{\operatorname{Dend}}
\newcommand{\prelie}{\operatorname{PreLie}}
\newcommand{\zinb}{\operatorname{Zinb}}
\newcommand{\leib}{\operatorname{Leib}}
\newcommand{\perm}{\operatorname{Perm}}

\newcommand{\ga}{\langle}
\newcommand{\dr}{\rangle}

\newtheorem{theorem}{Theorem}[section] 
 
\newtheorem{corollary}[theorem]{Corollary} 
\newtheorem{lemma}[theorem]{Lemma} 
\theoremstyle{definition}

\newtheorem{conjecture}[theorem]{Conjecture} 
 
\begin{document}
\title{On some anticyclic operads}
\author{F. Chapoton}
\address{Institut Girard Desargues, Universit\'e Claude Bernard 
(Lyon 1)\\B\^atiment Braconnier, 21 Avenue Claude Bernard\\F-69622 
Villeurbanne Cedex, France}

\asciiaddress{Institut Girard Desargues, Universite Claude Bernard 
(Lyon 1)\\Batiment Braconnier, 21 Avenue Claude Bernard\\F-69622 
Villeurbanne Cedex, France}

\email{chapoton@igd.univ-lyon1.fr}

\begin{abstract}
  Some binary quadratic operads are endowed with anticyclic structures
  and their characteristic functions as anticyclic operads are
  determined, or conjectured in one case.
\end{abstract}

\maketitle

\setcounter{section}{-1}

\section{Introduction}

The most classical operads are probably the three operads describing
commutative and associative algebras, associative algebras and Lie
algebras. They fit in the following diagram:
\begin{equation}
\xymatrix{
\comm & \asso \ar[l] & \lie \ar[l]
}
\end{equation}
These operads have very nice properties. First of all, they are binary
quadratic and Koszul in the sense of the work of Ginzburg and Kapranov
on the Koszul duality of operads \cite{ginzkapr}. More precisely
$\asso$ is self-dual and $\lie$ and $\comm$ are dual to each other.
Moreover they are basic examples of cyclic operads, a notion
introduced in \cite{getzkapr} and related to the moduli spaces of
curves with marked points. The developments of Koszul duality of
operads and of the theory of cyclic operads were both partly motivated
by the work of Kontsevich on non-commutative symplectic geometry
\cite{kont}, where three parallel constructions are made for these
three operads.

Our aim is to explain that most of the properties of this classical
sequence of operads also hold for two other diagrams involving some
binary quadratic operads introduced by Loday \cite{loday} and others
\cite{endo,rooted}. The first of these diagrams is
\begin{equation}
\xymatrix{
\perm & \dias \ar[l] & \leib \ar[l]
}
\end{equation}
and the other one is:
\begin{equation}
\xymatrix{
\zinb & \dend \ar[l] & \prelie \ar[l]
}
\end{equation}
All these operads are already known to be binary quadratic and Koszul.
More precisely, $(\perm,\prelie)$, $(\dias,\dend)$ and $(\leib,\zinb)$
are Koszul dual pairs of operads. The main objective of this article
is to show that they are anticyclic operads and that the maps in the
two diagrams above are maps of anticyclic operads.

One interesting point with cyclic and anticyclic operads is that they
allow the building of Lie algebras and graph complexes
\cite{modular,conant,markl}.  There has been a lot of work on the
graph complexes for the three classical cyclic operads $\comm$,
$\asso$ and $\lie$. It may be worth considering the analogous
structures for the six new anticyclic operads introduced here.

Let us also remark that there is a fourth classical cyclic operad,
namely the Poisson operad, which can be obtained as the graded cyclic
operad associated to a filtration of the cyclic operad $\asso$.
Similar objects exists in the diagrams above, namely filtrations of
the anticyclic operads $\dias$ and $\dend$ and associated graded
anticyclic operads, one of which is related to the pre-Poisson
algebras studied in \cite{aguiar}.

\section{Anticyclic operads}\label{general}

We briefly state some general facts on operads and anticyclic operads.
A convenient reference on this subject is \cite{markl}, see also
\cite{ginzkapr,getzkapr}. Most of the operads considered here will be
in the monoidal category of vector spaces over the field $\QQ$, but
the true ambient category is the category of chain complexes of vector
spaces over $\QQ$.

Recall that an operad $\PP$ is a collection of modules $\PP(n)$ over
the symmetric groups $\sym_n$ together with composition maps
satisfying some axioms modelled after the composition of multi-linear
maps. A non-symmetric operad $\PP$ is a similar structure without the
actions of the symmetric groups. If $\PP$ is a non-symmetric operad
then the collection $\PP(n)\otimes \QQ \sym_n$ is naturally an operad.

An anticyclic non-symmetric operad is a non-symmetric operad $\PP$
together with an action of a cyclic group of order $n+1$ on $\PP(n)$
satisfying some axioms. In particular, the action of the cyclic group
is determined by the action on the generators of the non-symmetric
operad. 

Similarly, an anticyclic operad is an operad $\PP$ together with an
action of a symmetric group $\sym_{n+1}$ on $\PP(n)$ extending the
action of $\sym_n$ and satisfying similar axioms. Here the group
$\sym_n$ is embedded in $\sym_{n+1}$ as the stabilizer of $n+1$. In
this case too, the action of the big symmetric group is determined by
the action on the generators of the operad.

The main difference with the notion of cyclic operad is that the unit
$1$ of the operad is mapped to $-1$ in an anticyclic case and to $1$
in the cyclic case. There is also some change of sign in the axioms.

The Hadamard product of two operads $\PP_1$ and $\PP_2$ is an operad
with underlying modules $\PP_1(n) \otimes \PP_2(n)$. If $\PP_1$ and
$\PP_2$ are cyclic or anticyclic operads, then the Hadamard product is
an anticyclic or cyclic operad for the tensor product action of the
symmetric groups $\sym_{n+1}$. Whether the Hadamard product is cyclic
or anticyclic is determined by the action on the unit $1$,
i.e.\ is given by a sign rule where ``cyclic'' is $+1$ and
``anticyclic'' is $-1$.

Let $\ZZ/{(n+1)\ZZ}$ be the subgroup of $\sym_{n+1}$ generated by the
longest cycle $(n+1,n,\dots,2,1)$. If $\PP$ is an anticyclic
non-symmetric operad, then the collection of induced modules
$\operatorname{Ind}_{\ZZ/{(n+1)\ZZ}}^{\sym_{n+1}}\PP(n)$ has a natural
structure of anticyclic operad.

There exists a differential graded anticyclic operad, called the
determinant operad and denoted by $\Det$, such that $\Det(n)$ is a
chain complex of dimension $1$ concentrated in degree $n-1$. The
Hadamard product by the operad $\Det$ corresponds to the suspension of
operads. It maps anticyclic operads to cyclic operads and vice-versa.

\section{Koszul duality and Legendre transform}

The theory of Koszul duality for operads has been introduced in
\cite{ginzkapr} for binary quadratic operads. 

For each binary quadratic operad $\PP$, one can define its dual operad
$\PP^!$ by elementary linear algebra using the quadratic presentation
of $\PP$. In particular, the space $\PP^!(2)$ of generators of $\PP^!$
is the tensor product of the dual space of $\PP(2)$ by the sign
representation of $\sym_2$.

A binary quadratic operad $\PP$ is Koszul if the natural morphisms of
cooperads from the bar complex of $\PP$ to the dual cooperad of the
suspension of $\PP^!$ is a quasi-isomorphism.

The Koszul dual operad $\PP^!$ of a Koszul operad $\PP$ has a natural
structure of cyclic (resp.\ anticyclic) operad if $\PP$ is anticyclic
(resp.\ anticyclic). When $\PP$ is cyclic or anticyclic, the action of
the cycle $(321)$ on $\PP^!(2)$ is given by the transpose of the
inverse of the action on the space $\PP(2)$.

When $\PP$ is cyclic (resp.\ anticyclic), both the bar complex and the
dual cooperad of the suspension of $\PP^!$ are anticyclic (resp.\
cyclic) and the natural map between them is automatically a morphism
of anticyclic (resp.\ cyclic) cooperads.

It is known that generating series of dual Koszul operads are related
by inversion for the plethysm. More precisely, let the characteristic
function of an operad $\PP$ be
\begin{equation}
  \ch(\PP)=\sum_{n \geq 1} \ch_n(\PP(n)),
\end{equation}
where $\ch_{n}(\PP(n))$ is the symmetric function for the
$\sym_{n}$-module $\PP(n)$.
The suspension $\Sigma F$ of a symmetric function $F$ is defined as
follows:
\begin{equation}
  \Sigma F=-F(-p_1,-p_2,-p_3,\dots),
\end{equation}
where the $p_k$ are power sums symmetric functions. Then one has the
following relation for a Koszul operad $\PP$:
\begin{equation}
  \ch(\PP) \circ \Sigma \ch(\PP^!)=p_1,
\end{equation}
where $\circ$ is the plethysm of symmetric functions.

A similar relation exists at the level of cyclic or anticyclic Koszul
operads, where the generating functions are related by the Legendre
transform of symmetric functions introduced in \cite{getzkapr}. Let us
give the precise statement.

The characteristic function of a cyclic or anticyclic operad $\PP$ is
defined by
\begin{equation}
   \Ch(\PP)=\sum_{n \geq 1} \Ch_{n+1}(\PP(n)),
\end{equation}
where $\Ch_{n+1}(\PP(n))$ is the symmetric function for the
$\sym_{n+1}$-module $\PP(n)$.

Let $F$ be a symmetric function with no term of degree $0$ and $1$ and
such that the term of degree $1$ of $\partial_{p_1}F$ does not vanish.
The Legendre transform $G=\legendre F$ of $F$ is defined by
\begin{equation}
  F \circ \frac{\partial G}{\partial p_1}+G=p_1 \frac{\partial G}{\partial p_1}.
\end{equation}
The Legendre transform is an involution, with the property that the
derived symmetric functions satisfy 
\begin{equation}
  \frac{\partial F}{\partial p_1} \circ \frac{\partial G}{\partial p_1}=p_1.
\end{equation}
Then if $\PP$ is an anticyclic (resp.\ cyclic) operad which is Koszul
as an operad, one has
\begin{equation}
  \legendre \Ch(\PP) = -\Sigma \Ch(\PP^!),
\end{equation}
where the anticyclic (resp.\ cyclic) structure on $\PP^!$ is the one
induced by the Koszul duality.

\section{The diassociative operad}

In this section and the next one, we will consider two non-symmetric
operads.

First let us recall some known facts about the non-symmetric operad
$\dias$ of diassociative algebras, see \cite{loday}.

First, it has a quadratic presentation. The generators are two binary
operations $x \dashv y$ and $x \vdash y$. These generators satisfy the
following relations:
\begin{eqnarray}
  x \dashv (y \vdash z)=x \dashv (y \dashv z) =(x \dashv y) \dashv z,\\
  \label{ax2} x \vdash (y \dashv z) =(x \vdash y) \dashv z,\\
  x \vdash (y \vdash z) =(x \vdash y) \vdash z=(x \dashv y) \vdash z.
\end{eqnarray}
Next, it is known that the space $\dias(n)$ has dimension $n$ with a
base \linebreak $(e^n_m)_{m=1,\dots,n}$ such that the composition in the 
operad is given by
\def\strut{\vrule width 0pt depth 8pt height 10pt}
\begin{equation}
  e^n_m \circ_i e^\ell_k =
  \begin{cases}
    \strut e^{n+\ell-1}_{m+k-1}&\text{ if }i=m,\\
    \strut e^{n+\ell-1}_{m+\ell-1}&\text{ if }i<m,\\
    \strut e^{n+\ell-1}_{m}&\text{ if }i>m.
  \end{cases}
\end{equation}
In the presentation above, $x_1 \dashv x_2$ is $e^2_1$ and $x_1 \vdash
x_2$ is $e^2_2$. More generally, the element $e^n_m$ is mapped to
\begin{equation}
  x_1 \vdash x_2 \vdash \dots \vdash x_m \dashv \dots \dashv x_n,
\end{equation}
with arbitrary parentheses. Conversely, any iterated product of the
variables $x_1,x_2,\dots,x_n$ in this order corresponds to some
$e^n_m$ by the following recursive procedure. At each step, choose the
sub-expression which is not on the bar side of $\dashv$ or $\vdash$,
until there remains only one variable $x_m$. For example,
\begin{equation}
  ((x_1 \dashv x_2) \vdash x_3) \dashv (x_4 \vdash x_5)
\end{equation}
is mapped to $e^5_3$.

Let us now introduce a notion of invariant bilinear map on a
diassociative algebra. It is an antisymmetric map with value in some
vector space:
\begin{equation}
  \ga x\, ,\, y \dr =-\ga y\, ,\, x \dr,
\end{equation}
such that
\begin{equation}
  \ga x \dashv y \, ,\, z \dr = -\ga y \vdash z \, ,\, x \dr
\text{  and  }
  \ga x \vdash y \, ,\, z \dr =
  \ga y \dashv z \, ,\, x \dr -\ga y \vdash z \, ,\, x \dr.
\end{equation}
Let us define a map $\tau_1$ as multiplication by $-1$ on $\dias(1)$
and a map $\tau_2$ on $\dias(2)$ by the following relation
\begin{equation}
  \ga E(x,y)\, ,\,z\dr = \ga \tau_2(E)(y,z)\, ,\,x \dr
\end{equation}
for each element $E$ in $\{e^2_1,e^2_2\}$, using the previous conditions
on $\ga \,,\,\dr$.

\begin{theorem}
  There exists a unique collection $\tau$ of endomorphisms $\tau_n$ of
  the space $\dias(n)$ extending $\tau_1$ and $\tau_2$ and endowing
  the operad $\dias$ with a structure of non-symmetric anticyclic
  operad.
\end{theorem}

\begin{proof}
  Clearly $\tau_1$ is of order $2$ and $\tau_2$ is of order $3$. As
  $\dias$ is generated by $\dias(2)$, the structure map $\tau$ is
  unique if it exists. To check that $\tau_n$ can be defined for
  $n\geq 3$, it is enough to check that the notion of invariant form
  is compatible with the relations defining the operad $\dias$. Let us
  check one of these compatibility conditions, for the relation
  (\ref{ax2}):
  \begin{eqnarray}
    & \ga x \vdash (y \dashv z)- (x \vdash y) \dashv z \,,\, t \dr \\ 
    &=
    \ga (y \dashv z) \dashv t- (y \dashv z) \vdash t \,,\, x \dr
     + \ga z \vdash t \,,\, x\vdash y \dr\\
    &=
    \ga (y \dashv z) \dashv t- (y \dashv z) \vdash t \,,\, x \dr
     - \ga x \vdash y \,,\, z\vdash t \dr\\
    &=  \ga (y \dashv z) \dashv t- (y \dashv z) \vdash t \,,\, x \dr
    - \ga y \dashv (z \vdash t)- y \vdash (z \vdash t) \,,\, x \dr =0,
  \end{eqnarray}
  where one has used the antisymmetry and the invariance to obtain an
  expression with the $x$ variable only in the right slot.
  The remaining checks are quite similar and are left to the reader.
\end{proof}

In fact, $\tau_n$ has a very simple expression in the base $e^n$.

\begin{lemma}
  The action of $\tau_n$ is given by
  \begin{equation}
    \label{companion}
    \tau_n(e^{n}_{m})=
    \begin{cases}
      -e^{n}_{n} &\text{ if }m=1,\\
      -e^{n}_{n}+e^{n}_{m-1} &\text{ else}.
    \end{cases}
  \end{equation}
\end{lemma}

\begin{proof}
  It is readily true for $n=1$ and $2$. Let us first prove that
  $\tau_{n+1}(e^{n+1}_{1})=e^{n+1}_{n+1}$. This follows from the
  equality:
  \begin{align}
    \ga x_1 \dashv x_2 \dashv \dots \dashv x_{n+1} \,,\,x_{n+2} \dr
    &=-\ga (x_2 \dashv \dots \dashv x_{n+1}) \vdash x_{n+2} \,,\, x_1 \dr\\
    &=-\ga x_2 \vdash \dots \vdash x_{n+2} \,,\, x_1 \dr.
  \end{align}
  Let us now compute $\tau_{n+1}(e^{n+1}_{m})$ for $m \geq 2$. Indeed,
  one has
  \begin{multline}
    \ga x_1 \vdash x_2 \vdash \dots \vdash x_m \dashv \dots \dashv
    x_{n+1} \,,\,x_{n+2} \dr=\\
    \ga ( x_2 \vdash \dots \vdash x_m \dashv \dots \dashv
    x_{n+1}) \dashv x_{n+2}
    - ( x_2 \vdash \dots \vdash x_m \dashv \dots \dashv
    x_{n+1}) \vdash x_{n+2} \,,\, x_1 \dr\\
    = \ga x_2 \vdash \dots \vdash x_m \dashv \dots \dashv x_{n+2}
    -   x_2 \vdash \dots \vdash x_{n+2} \,,\, x_1 \dr.\\
  \end{multline}
  This concludes the proof.
\end{proof}

One can note that the matrix of $\tau_n$ in the base $e^n$ is a
companion matrix for the polynomial $1+x+\dots+x^{n}$ and can also be
described as $-({}^{t}L)^{-1} L$ where $L$ is the lower triangular
matrix with $1$ everywhere below the diagonal. 

\section{The dendriform operad}

Let us now recall some facts about the non-symmetric operad $\dend$
describing dendriform algebras, see \cite{loday}.

First, it also has a quadratic presentation. The generators are two
binary operations $x \prec y$ and $x \succ y$. These generators
satisfy the following relations:
\begin{eqnarray}
  (x \prec y) \prec z)=x \prec (y \prec z)+x \prec (y \succ z),\\
  x \succ (y \prec z) =(x \succ y) \prec z,\\
  x \succ (y \succ z)=(x \succ y) \succ z+(x \prec y) \succ z.
\end{eqnarray}
Next, it is known that the dimension of the space $\dend(n)$ is the
Catalan number $c_n=\frac{1}{n+1}\binom{2n}{n}$. There is a base
$Y^n$ of $\dend(n)$ indexed by planar binary trees with $n+1$ leaves,
in which the composition of the operad has a simple expression. The
base $Y^2$ is made precisely of $x_1 \prec x_2$ and $x_1 \succ x_2$.

Let us introduce a notion of invariant bilinear map on a dendriform
algebra. It is an antisymmetric map:
\begin{equation}
  \ga x\, ,\, y \dr =-\ga y\, ,\, x \dr
\end{equation}
such that
\begin{equation}
  \ga x \succ y \, ,\, z \dr =\ga y \prec z \, ,\, x \dr
  \text{  and  }
  \ga x \prec y \, ,\, z \dr = -\ga y \prec z \, ,\, x
  \dr -\ga y \succ z \, ,\, x
  \dr.
\end{equation}
As before, let us define a map $\tau_1$ as multiplication by $-1$ on
$\dend(1)$ and a map $\tau_2$ on $\dend(2)$ by the following relation
\begin{equation}
  \ga E(x,y)\, ,\,z\dr = \ga \tau_2(E)(y,z)\, ,\,x \dr
\end{equation}
for each element $E$ in the base $Y^2$ of $\dend(2)$, using the
previous conditions on $\ga \,,\,\dr$.

\begin{theorem}
  There exists a unique collection $\tau$ of endomorphisms $\tau_n$ of
  the space $\dend(n)$ extending $\tau_1$ and $\tau_2$ and endowing
  the operad $\dend$ with a structure of non-symmetric anticyclic
  operad.
\end{theorem}

\begin{proof}
  Clearly $\tau_1$ is of order $2$ and $\tau_2$ is of order $3$. As
  $\dend$ is generated by $\dend(2)$, the maps $\tau$ are unique if
  they exist. To check that $\tau_n$ can be defined for $n\geq 3$, it
  is enough to check that the notion of invariant form is compatible
  with the three relations defining the operad $\dend$.  Let us check
  that it works for the second compatibility condition:
  \begin{eqnarray}
    & \ga x \succ (y \prec z)- (x \succ y) \prec z \,,\, t \dr \\ 
    &=
    \ga (y \prec z) \prec t \,,\, x \dr
     + \ga z \prec t \,,\, x\succ y \dr + \ga z \succ t \,,\, x\succ y \dr\\
    &=
    \ga (y \prec z) \prec t \,,\, x \dr
     - \ga x\succ y \,,\, z \prec t \dr - \ga x\succ y\,,\, z \succ t  \dr\\
    &=  \ga (y \prec z) \prec t \,,\, x \dr
    - \ga y \prec (z \prec t)+ y \prec (z \succ t) \,,\, x \dr =0,
  \end{eqnarray}
  where one has used the antisymmetry and the invariance to obtain an
  expression with the $x$ variable alone on the right.  The two
  remaining computations are just as simple and are left to the
  reader.
\end{proof}

The matrix of $\tau$ in the base of trees certainly deserves more
study. It seems to be related to the so-called Tamari lattices
\cite{tamari}.

\section{Four operads}

Starting from this section, we consider operads in the usual sense,
which means with actions of the symmetric groups. As explained in
section \ref{general}, we can consider the two non-symmetric
anticyclic operads just defined as anticyclic operads, still denoted
$\dias$ and $\dend$. We will show that some sub-operads and quotient
operads of these inherit an anticyclic structure.

\subsection{The Perm operad}

The $\perm$ operad, introduced in \cite{endo}, is a quotient operad of
the diassociative operad $\dias$. The space $\perm(n)$ has dimension
$n$ and the action of $\sym_n$ is the usual permutation
representation.

The operad $\perm$ is the quotient of $\dias$ by the ideal generated
by the element $x_1 \dashv x_2 - x_2 \vdash x_1 $. The image of the
product $x_1 \dashv x_2$ will be denoted $x_1 x_2$.

The operad $\perm$ is quadratic, generated by the binary product $xy$
(regular representation of $\sym_2$) and with relations
\begin{equation}
  (xy)z=x(yz)=x(zy).
\end{equation}

\begin{theorem}
  The operad $\perm$ has a unique structure of anticyclic operad such
  that the quotient map from $\dias$ is a morphism of anticyclic
  operads.
\end{theorem}

\begin{proof}
  One has to check that the defining ideal is stable by the action of
  $\sym_{n+1}$ on $\dias(n)$. It is enough to check this in
  $\dias(2)$, which contains the generators of the ideal, where it is
  immediate.  Hence $\perm$ is a quotient anticyclic operad of
  $\dias$.
\end{proof}

The resulting notion of invariant bilinear map is as follows. It is an
antisymmetric map
\begin{equation}
  \ga x\, ,\, y \dr =-\ga y\, ,\, x \dr
\end{equation}
such that
\begin{equation}
   \ga x y \, ,\, z \dr = -\ga z y\, ,\, x  \dr
\text{  and  }
    \ga x  y \, ,\, z \dr 
  = \ga x  z \, ,\, y \dr - \ga z x \, ,\, y \dr.
\end{equation}

\begin{theorem}
  \label{reflex}
  The action of $\sym_{n+1}$ on $\perm(n)$ is isomorphic to the
  representation by reflections.
\end{theorem}

\begin{proof}
  Let us consider the action of $\sym_{n+1}$ by permutations on the
  module with base $(\varepsilon_i)_{i=1,\dots,n+1}$. The reflection
  module is the submodule with base
  $b^n_m=\varepsilon_m-\varepsilon_{n+1}$ for $m=1,\dots,n$. The
  action of the subgroup $\sym_n$ is by permutations of the vectors
  $b^n_m$ and the action of the cycle $\tau_{n}$: $(n+1,n,\dots,2,1)$
  is given by
   \begin{equation}
    \tau_n(b^{n}_{m})=
    \begin{cases}
      -b^{n}_{n} &\text{ if }m=1,\\
      -b^{n}_{n}+e^{n}_{m-1} &\text{ else}.
    \end{cases}
  \end{equation}
  On the other hand, the action of $\sym_n$ on the module $\perm(n)$
  with base \linebreak $(e^n_m)_{m=1,\dots,n}$ is by permutation of the 
  vectors
  $e^n_m$. This base of $\perm(n)$ is the image of the base $e^n_m
  \otimes 1$ of the operad $\dend$. Hence the action of the cycle
  $\tau_n$ is induced by the action of the cycle $\tau_n$ in the base
  $e^n_m$ of the non-symmetric operad $\dias$ which is given by 
  Formula (\ref{companion}). This concludes the proof.
\end{proof}

\begin{corollary}
  The characteristic function of the anticyclic operad $\perm$ is
  \begin{equation}
    \label{ch_perm}
    \Ch(\perm)=(p_1-1) \exp( \sum_{k \geq 1} \frac{p_k}{k} ) +1.
  \end{equation}
\end{corollary}
\begin{proof}
  This follows readily from the previous Theorem and the well-known
  fact that
  \begin{equation}
    \ch(\comm)=\exp( \sum_{k \geq 1} \frac{p_k}{k} )-1.
  \end{equation}
\vspace{-25pt}

\end{proof}

\subsection{The Leibniz operad}

The Leibniz operad is a sub-operad of $\dias$, also introduced in
\cite{loday}. It will be denoted by $\leib$.

It is the sub-operad of $\dias$ generated by the element
$[x_1,x_2]=x_1 \dashv x_2 - x_2 \vdash x_1 $. Beware that this bracket
is not antisymmetric. Leibniz algebras are ``non-commutative Lie
algebras'' in some sense.

The quadratic presentation of $\leib$ is the following. It is
generated by the binary product $[x,y]$ (regular module of $\sym_2$)
modulo the relation
\begin{equation}
  [x,[y,z]]=[[x,y],z]-[[x,z],y].
\end{equation}
It is known that the space $\leib(n)$ is the regular representation of
$\sym_n$, with base given by left-bracketed words in the variables
$x_1,\dots,x_n$.

\begin{theorem}
  The operad $\leib$ has a unique structure of anticyclic operad
  such that the inclusion map into $\dias$ is a morphism of anticyclic
  operads.
\end{theorem}

\begin{proof}
  It is enough to check that $\leib(2)$ is indeed a submodule of
  $\dias(2)$ for the action of $\sym_3$. This is immediate and implies
  that $\leib(n)$ is stable by the action of the symmetric group
  $\sym_{n+1}$ on $\dias(n)$.
\end{proof}

The resulting notion of invariant bilinear map is as follows. It is an
antisymmetric map
\begin{equation}
  \ga x\, ,\, y \dr =-\ga y\, ,\, x \dr
\end{equation}
such that
\begin{equation}
  \ga [x,  y] \, ,\, z \dr = \ga [z ,y] \, ,\, x  \dr
\text{  and  }
  \ga [x, y] \, ,\, z \dr 
  =- \ga [z,x] \, ,\, y \dr - \ga [x,z] \, ,\, y \dr.
\end{equation}

\begin{theorem}
  \label{ch_leib}
  The action of $\sym_{n+1}$ on $\leib(n)$ is isomorphic to the Lie
  module $\lie(n+1)$.
\end{theorem}

\begin{proof}
  It was proved in \cite{endo} that the operad $\leib$ is the Hadamard
  product of the operads $\perm$ and $\lie$. As explained in Section
  \ref{general}, the Hadamard product of a cyclic operad and an anticyclic
  operad is an anticyclic operad. One can check that the anticyclic
  structure obtained in this way on $\leib$ coincides with the one
  introduced above. Hence the action of $\sym_{n+1}$ on $\leib(n)$ is
  given by the tensor product of the $\sym_{n+1}$-modules $\lie(n)$
  and $\perm(n)$. Then using Theorem \ref{reflex} and \cite[Corollary
  6.8]{getzkapr}, this is known to be isomorphic to the module
  $\lie(n+1)$.
\end{proof}

\subsection{The Pre-Lie operad}

The $\prelie$ operad has been introduced in \cite{rooted}. It is the
sub-operad of the dendriform operad generated by the operation $x_1
\sous x_2 = x_2 \succ x_1-x_1 \prec x_2$. The space $\prelie(n)$ has
dimension $n^{n-1}$.

Let us recall the presentation of the operad $\prelie$. The product $x
\sous y$ is the regular module for $\sym_2$ and must satisfy the
following relation:
\begin{equation}
  (x\sous y)\sous z- x \sous (y \sous z)
  =(x\sous z)\sous y- x \sous (z \sous y).
\end{equation}

\begin{theorem}
  The operad $\prelie$ has a unique structure of anticyclic operad
  such that the inclusion map into $\dend$ is a morphism of anticyclic
  operads.
\end{theorem}

\begin{proof}
  It is enough to check that $\prelie(2)$ is indeed a submodule of
  $\dend(2)$ for the $\sym_3$ action. This is immediate and implies
  that $\prelie(n)$ is stable by the action of the symmetric group
  $\sym_{n+1}$ on $\dend(n)$.
\end{proof}

The resulting notion of invariant bilinear map is as follows. It is an
antisymmetric map
\begin{equation}
  \ga x\, ,\, y \dr =-\ga y\, ,\, x \dr
\end{equation}
such that
\begin{equation}
  \ga x \sous y \, ,\, z \dr = -\ga x \sous  z \, ,\, y
  \dr
\text{  and  }
  \ga x \sous y \, ,\, z \dr 
  =- \ga y \sous z \, ,\, x \dr+ \ga z \sous y \, ,\, x \dr.
\end{equation}

{\bf Remark}\qua As the space $\prelie(n)$ is isomorphic as a $\sym_n$-module
to the space with a base indexed by labelled rooted trees on $n$
vertices, this implies the existence of a remarkable linear action of
$\sym_{n+1}$ on this space.

\begin{theorem}
  One has the following equality of symmetric functions:
  \begin{equation}
    \label{bizarre}
    \ch(\prelie)(1+\Ch(\prelie))=p_1(1+\ch(\prelie)+\ch(\prelie)^2). 
  \end{equation}
\end{theorem}
\begin{proof}
  One can check that $\prelie$ is the Koszul dual anticyclic operad of
  $\perm$. Hence its characteristic function is obtained by a Legendre
  transform of the characteristic function of $\perm$.

  But since $\perm$ and $\prelie$ are Koszul dual operads, it is known
  that
  \begin{equation}
    \ch(\perm) \circ \Sigma \ch(\prelie)=p_1.
  \end{equation}
  This can be used to replace the equation defining the Legendre
  transform of $\Ch(\perm)$ by a relation no longer involving the
  plethysm. Applying the suspension $\Sigma$ to this relation leads to
  Formula (\ref{bizarre}).
\end{proof}

We propose here an explicit conjecture for the character of
$\sym_{n+1}$ on $\prelie(n)$ as a symmetric function.

\begin{conjecture}
  \label{conj_prelie}
  The characteristic function $\Ch(\prelie)$ of the anticyclic operad
  $\prelie$ is
  \begin{equation}
    \sum_{\lambda,|\lambda|\geq 1,\lambda_1\not=1}(\lambda_1-1)^{\lambda_1-2}\prod_{k \geq
      2}\left( (f_k(\lambda)-1)^{\lambda_k}-k \lambda_k (f_k(\lambda)-1)^{\lambda_k-1}\right) \frac{p_\lambda}{z_\lambda},
  \end{equation}
  where the sums runs over non-empty partitions $\lambda$ not having
  exactly one part of size $1$, $\lambda_k$ denotes the number of
  parts of size $k$ in the partition $\lambda$ and $f_k(\lambda)$
  denotes the number of fixed points of the $k^{th}$ power of a
  permutation of cycle type $\lambda$. The notations $p_{\lambda}$ and
  $z_{\lambda}$ are classical for power sum symmetric functions and
  related constants.
\end{conjecture}

It is easy to check that the restriction to $\sym_n$ gives back the
known formula for the action on rooted trees which can be found in
\cite{labelle}. It has been checked by computer up to $n=14$ that the
expected characteristic function is indeed a positive linear
combination of Schur functions and that Formula (\ref{bizarre}) holds.

\subsection{The Zinbiel operad}

The Zinbiel operad, denoted by $\zinb$, was introduced in \cite{loday}.
Maybe it would be more appealing to call it the shuffle operad. It is
the quotient operad of $\dend$ by the ideal generated by the following
element
\begin{equation}
  x_1 \prec x_2 - x_2 \succ x_1.
\end{equation}
The image in $\zinb$ of the product $x_1 \prec x_2$ will be denoted
$x_1 x_2$.

The operad $\zinb$ has a quadratic presentation. It is generated by
the binary product $xy$ (regular representation of $\sym_2$) subject
to the relation:
\begin{equation}
  (xy)z=x(yz)+x(zy).
\end{equation}
The space $\zinb(n)$ is isomorphic to the regular representation of
$\sym_n$ and the composition of the operad can be described using
shuffles of permutations.

\begin{theorem}
  The operad $\zinb$ has a unique structure of anticyclic operad such
  that the quotient map from $\dend$ is a morphism of anticyclic
  operads.
\end{theorem}

\begin{proof}
  Once again, it follows already from the check done for the $\prelie$
  operad that the ideal defining $\zinb$ is indeed stable by the
  action of the symmetric group $\sym_{n+1}$ on $\dend(n)$. Hence
  $\zinb$ is a quotient anticyclic operad of $\dend$.
\end{proof}

The resulting notion of invariant bilinear map is as follows. It is an
antisymmetric map
\begin{equation}
  \ga x\, ,\, y \dr =-\ga y\, ,\, x \dr
\end{equation}
such that
\begin{equation}
  \ga x y \, ,\, z \dr = \ga x z \, ,\, y  \dr
\text{  and  }
  \ga x  y \, ,\, z \dr 
  =- \ga y  z \, ,\, x \dr - \ga z y \, ,\, x \dr.
\end{equation}

\begin{theorem}
  The action of $\sym_{n+1}$ on $\zinb(n)$ is isomorphic to the Lie
  module $\lie(n+1)$.
\end{theorem}

\begin{proof}
  It is known that the operad $\zinb$ is the Koszul dual of the operad
  $\leib$. One can check that the anticyclic structure of $\leib$ is
  obtained from the anticyclic structure of $\zinb$ by Koszul duality.
  Hence it follows that the characteristic of the anticyclic operad
  $\zinb$ is related to the characteristic of the anticyclic operad
  $\leib$ by a Legendre transform of symmetric functions. The
  characteristic of $\leib$ is known by Theorem \ref{ch_leib} to be
  \begin{equation}
    F=\sum_{n \geq 2} \frac{1}{n} \sum_{d | n} \mu(d) p_d^{n/d}. 
  \end{equation}
  One has to check that the Legendre transform of $F$ is
  \begin{equation}
    G=\sum_{n \geq 2} \frac{1}{n} \sum_{d | n} \mu(d) (-1)^{n/d} p_d^{n/d}. 
  \end{equation}
  The derivative of $G$ is
  \begin{equation}
    \partial_{p_1}G=\frac{p_1}{1+p_1},
  \end{equation}
  and one has
  \begin{equation}
    p_1\partial_{p_1}G=\frac{p_1^2}{1+p_1}.
  \end{equation}
  Let us compute $F \circ \partial_{p_1}G +G$. One finds
  \begin{equation}
     \sum_{n \geq 2} \frac{1}{n} \sum_{d | n} \mu(d) 
     \left(\frac{1}{(1+p_d)^{n/d}}+(-1)^{n/d}\right) p_d^{n/d}. 
  \end{equation}
  Now fix $d \geq 2$ and look only at the extracted series in $p_d$:
  \begin{equation}
    \frac{\mu(d)}{d} \sum_{k \geq 1} \frac{1}{k} \left( \frac{1}
      {(1+p_d)^k}+(-1)^k \right) p_d^k.
  \end{equation}
  This is easily expressed  using logarithms and seen to vanish.
  Hence the expression $F \circ \partial_{p_1}G +G$ is given by
  \begin{equation}
     \sum_{n \geq 2} \frac{1}{n}  
     \left(\frac{1}{(1+p_1)^{n}}+(-1)^{n}\right) p_1^{n}=\frac{p_1^2}{1+p_1} 
   \end{equation}
  as expected.
\end{proof}

Let us remark that exactly the same computation with the M\"obius
function $\mu$ replaced by the Euler totient function $\varphi$
corresponds to self-duality of the characteristic function for the
cyclic Koszul operad $\asso$.

\section{Characters of $\rm Dend$ and $\rm Dias$}

\begin{theorem}
  The characteristic function of the anticyclic operad $\dias$ is
  \begin{equation}
    \label{ch_dias}
    \Ch(\dias)=
    \sum_{n \geq 2} \left(p_1^n-\frac{1}{n}\sum_{d | n} \varphi(d) p_{d}^{n/d}\right).
  \end{equation}
\end{theorem}

\begin{proof}
  It is easy to check that the anticyclic operad $\dias$ is obtained
  as the Hadamard product of the anticyclic operad $\perm$ by the
  cyclic operad $\asso$.

  It is known that the characteristic function of the cyclic operad
  $\asso$ is
  \begin{equation}
    \sum_{n \geq 2}\frac{1}{n}\sum_{d | n} \varphi(d) p_{d}^{n/d},
  \end{equation}
  where $\varphi$ is the Euler totient function.
  
  Using Theorem \ref{reflex}, one knows that the tensor product by the
  $\sym_{n+1}$ module $\perm(n)$ is given by the operator
  $-\id+\partial_{p_1}$. Then a simple computation proves the Theorem.
\end{proof}

One can check that the anticyclic structure of $\dend$ is the one
obtained from the anticyclic structure of $\dias$ by Koszul duality.
Hence, by Legendre inversion, one gets that the characteristic
function of $\dend$ is related to the Legendre transform of the
characteristic function of $\dias$.

\begin{theorem}
  The characteristic function of the anticyclic operad $\dend$ is 
  \begin{equation}
    \Ch(\dend)=1-p_1-\sqrt{1-4 p_1} -\sum_{n \geq 1} \left(\frac{1}{2n} \sum_{d | n} 
    \varphi(d) \binom{2n/d}{n/d}p_{d}^{n/d}\right).
  \end{equation}
\end{theorem}
\begin{proof}
  Let us check that the Legendre transform gives this result. Let $F$
  be the characteristic function (\ref{ch_dias}) of $\dias$ and let
  $G$ be
  \begin{equation}
    1+p_1-\sqrt{1+4 p_1} -\sum_{n \geq 1} \frac{1}{2n} \sum_{d | n} 
    \varphi(d) \binom{2n/d}{n/d}(-p_{d})^{n/d}.
  \end{equation}
  One has to check that $F \circ \partial_{p_1} G+G=p_1\partial_{p_1} G$.

  Then one has
  \begin{equation}
    \partial_{p_1}G=\frac{1+2 p_1-\sqrt{1+4 p_1}}{2 p_1}
  \end{equation}
  and
  \begin{equation}
    p_1\partial_{p_1}G=\frac{1+2 p_1-\sqrt{1+4 p_1}}{2}.
  \end{equation}
  Let us compute $p_1 \partial_{p_1}G-F \circ \partial_{p_1}G$. One finds
  \begin{multline}
    \frac{1+2 p_1-\sqrt{1+4 p_1}}{2}-\left(\frac{1+2 p_1-\sqrt{1+4
          p_1}}{2 p_1}\right)^2\frac{1}{1-\left(\frac{1+2
          p_1-\sqrt{1+4 p_1}}{2 p_1}\right)} \\+\sum_{n \geq 2}
    \frac{1}{n} \sum_{d | n} \varphi(d) \left(\frac{1+2 p_d-\sqrt{1+4
          p_d}}{2 p_d}\right)^{n/d}.
  \end{multline}
  Let us split this sum into
  \begin{equation}
    \frac{1+2 p_1-\sqrt{1+4 p_1}}{2}-\frac{\left(1+2 p_1-\sqrt{1+4
          p_1}\right)^2}{2 p_1 \left(\sqrt{1+4
          p_1}-1\right)}-\frac{1+2 p_1-\sqrt{1+4 p_1}}{2 p_1}
  \end{equation}
  and
  \begin{equation}
    \sum_{n \geq 1} \frac{1}{n} \sum_{d | n} 
    \varphi(d) \left(\frac{1+2 p_d-\sqrt{1+4 p_d}}{2 p_d}\right)^{n/d}.
  \end{equation}
  Let us compute these separately. The first part is easily seen to be
  \begin{equation}
    1+ p_1-\sqrt{1+4 p_1}.
  \end{equation}
  The second part becomes
  \begin{equation}
    \sum_{d \geq 1} \frac{\varphi(d)}{d} \sum_{k \geq 1} \frac{1}{k}
    \left(\frac{1+2 p_d-\sqrt{1+4 p_d}}{2 p_d} \right)^k,
  \end{equation}
  which is
  \begin{equation}
    \sum_{d \geq 1} \frac{\varphi(d)}{d} 
    \log\left(\frac{1+\sqrt{1-4p_d}}{2}\right).
  \end{equation}
  Using then the Taylor expansion
  \begin{equation}
    -\log\left(\frac{1+\sqrt{1-4u}}{2}\right)
    =\sum_{n\geq 1} \frac{1}{2n}\binom{2n}{n}u^n,
  \end{equation}
  one gets
  \begin{equation}
    -\sum_{d \geq 1} \frac{\varphi(d)}{d} 
    \sum_{k \geq 1} \frac{1}{2k}\binom{2k}{k}(-p_d)^k,
  \end{equation}
  which becomes
  \begin{equation}
    -\sum_{n \geq 1} \frac{1}{2n} \sum_{d |n} \varphi(d) 
    \binom{2n/d}{n/d}(-p_d)^{n/d}.
  \end{equation}
  Summing both parts gives, as expected, that $p_1\partial_{p_1} G-F
  \circ \partial_{p_1} G=G$.
\end{proof}

\Addressesr

\end{document}